\documentclass[12pt,letterpaper]{amsart}
\usepackage{amsmath,amsfonts,amssymb,amsthm,amscd}
\usepackage[mathscr]{eucal}
\usepackage{graphicx}
\usepackage{appendix}
\usepackage[utf8]{inputenc}
\usepackage[T1]{fontenc}
\usepackage{lmodern}
\usepackage[brazil]{babel}

\newcommand{\e}{\epsilon}

\newcommand{\real}{\mathbb{R}}

\newcommand{\torus}{{\mathbb T}}

\newcommand{\dd}{\,\mathrm{d}}

\DeclareMathOperator{\dv}{div} %
\DeclareMathOperator{\rote}{rot}%

\newcommand{\loc}{\operatorname{{loc}}}

\newtheorem{theorem}{Theorem}

\theoremstyle{definition}
\newtheorem{definition}{Definition}

\theoremstyle{remark}

\def \e1 {\vec{e}_1}
\def \e2 {\vec{e}_2}

\begin{document}


\title[Navier-Stokes e o Problema do Mil\^{e}nio]{Singularidades para as solu\c{c}\~{o}es das Equa\c{c}\~{o}es de Navier-Stokes and Euler e o Problema do Mil\^{e}nio}

\author[M. C. Lopes Filho]{Milton C. Lopes Filho}
\address[Milton C. Lopes Filho]{Instituto de Matem\'atica\\Universidade Federal do Rio de Janeiro\\Cidade Universit\'aria -- Ilha do Fund\~ao\\Caixa Postal 68530\\21941-909 Rio de Janeiro, RJ -- BRAZIL.}
\email{mlopes@im.ufrj.br}\pagestyle{myheadings}
\author[H. J. Nussenzveig Lopes]{Helena J. Nussenzveig Lopes}
\address[Helena J. Nussenzveig Lopes]{Instituto de Matem\'atica\\Universidade Federal do Rio de Janeiro\\Cidade Universit\'aria -- Ilha do Fund\~ao\\Caixa Postal 68530\\21941-909 Rio de Janeiro, RJ -- BRAZIL.}
\email{hlopes@im.ufrj.br}

\date{\today}

\begin{abstract}
O objetivo desta nota \'{e} oferecer um panorama da hist\'{o}ria e do estado da arte na pesquisa em torno do Problema do Mil\^{e}nio para as equa\c{c}\~{o}es de Navier-Stokes incompressiveis, singularidades em din\^{a}mica dos fluidos, e o problema correspondente para as equa\c{c}\~{o}es de Euler.
\end{abstract}

\maketitle

Escoamento \'{e} um termo que se refere \`{a} descri\c{c}\~{a}o de fluidos em movimento. Fluidos, e escoamentos, aparecem no mundo real de in\'{u}meras maneiras, desde escoamentos intracelulares, em capilares ou meios porosos, at\'{e} movimentos em escalas planet\'{a}rias, tais como a atmosfera de Júpiter.  Escoamentos também apresentam diversas características físicas distintas, dependendo da maneira com que os fluidos respondem a for\c{c}as aplicadas. Por exemplo, em escoamentos incompressíveis, o volume ocupado por uma parcela de fluido varia muito pouco ao ser comprimido. Em escoamentos Newtonianos, os fluidos s\~{a}o indiferentes a deforma\c{c}\~{a}o, resistindo apenas \`{a} velocidade com que deforma\c{c}\~{o}es s\~{a}o impostas. As equa\c{c}\~{o}es de Navier-Stokes s\~{a}o o modelo matem\'{a}tico fundamental para o movimento de fluidos Newtonianos  e incompressíveis, servindo bem para estudar escoamentos em uma enorme gama de situa\c{c}\~{o}es de interesse te\'{o}rico e pr\'{a}tico.

As equa\c{c}\~{o}es de Navier-Stokes para um fluido incompressível e Newtoniano, que ocupa o espa\c{c}o $N$-dimensional todo, se escrevem, de forma condensada como:

\begin{equation} \label{EDPsFluidos}
\left\{
\begin{array}{ll}
\partial_t u + (u \cdot \nabla) u = - \nabla p +  \nu \Delta u + F, & \mbox{ em } \mathbb{R}^N \times (0,\infty)\\
\dv u = 0, & \mbox{ em } \mathbb{R}^N \times [0,\infty)\\
|u(x,t)| \to 0, & \mbox{ quando } |x| \to \infty, \; t>0\\
u(x,0)=u_0(x), & \mbox{ em } \mathbb{R}^N \times \{t=0\}.
\end{array}\right.
\end{equation}

\vspace{0.5cm}

Acima, o campo  $u = (u_1, u_2, \ldots, u_N)$ representa a velocidade do escoamento no instante $t$ e no ponto $x \in \real^N$; $p=p(x,t)$ \'{e} \ a press\~{a}o escalar. O termo $(u \cdot \nabla) u$ \'{e} \ entendido como sendo:
\[[(u \cdot \nabla) u]_j = \sum_{i=1}^N u_i \partial_{x_i} u_j,\;\;\;j=1,\ldots, N.\]
O termo $F=F(x,t)$ \'{e} uma for\c{c}a externa, por exemplo, $F=-g e_N$ representa a for\c{c}a devida \`{a} gravidade. O parâmetro $\nu > 0$ \'{e} a viscosidade do fluido - uma medida da fric\c{c}\~{a}o interna do fluido.

A primeira equa\c{c}\~{a}o em \eqref{EDPsFluidos} (na verdade $N$ equa\c{c}\~{o}es escritas vetorialmente) exprime a {\it segunda lei de Newton}, em uma situa\c{c}\~{a}o em que a densidade do fluido é constante e normalizada com o valor unit\'{a}rio. A segunda equa\c{c}\~{a}o exprime a incompressibilidade do escoamento. O sistema tem $N+1$ equa\c{c}\~{o}es e $N+1$ inc\'{o}gnitas, a saber, as $N$ componentes da velocidade e a press\~{a}o escalar $p$.

No caso particular $\nu = 0$, o sistema \eqref{EDPsFluidos} torna-se um modelo de {\it fluidos ideais ou invíscidos}, e o sistema \'{e} chamado {\it equa\c{c}\~{o}es de Euler}. Sem fric\c{c}\~{a}o, um fluido ideal, uma vez posto em movimento, nunca mais para.

Se $\nu > 0$, as solu\c{c}\~{o}es de \eqref{EDPsFluidos} representam escoamentos viscosos, e a fric\c{c}\~{a}o tende a fazer o movimento parar, dissipando a energia presente inicialmente. Todo escoamento no mundo real tem alguma viscosidade.

No que se segue, vamos nos concentrar nos casos de interesse físico, para os quais $N=2,3$. Muito de nossa discussão é baseada no excelente texto de G. Lemari\'{e}-Rieusset \cite{LM-Rbook}, que convidamos o leitor interessado a consultar.

As equa\c{c}\~{o}es de Euler apareceram pela primeira vez no artigo {\it Principes g\'{e}n\'{e}raux du mouvement des fluides}, publicado por L. Euler em $1757$, \cite{Euler1757}. Foi a segunda equa\c{c}\~{a}o diferencial parcial a aparecer em uma publica\c{c}\~{a}o. A primeira foi a equa\c{c}\~{a}o da onda, formulada por J. D'Alembert em $1747$, \cite{DAlembert}. As equa\c{c}\~{o}es de Navier-Stokes tiveram uma hist\'{o}ria mais complicada, aparecendo pela primeira vez em um artigo de Navier em $1822$, \cite{Navier} sendo depois redescobertas, ou reinventadas com diferentes justificativas, por Cauchy em $1823$, \cite{Cauchy1,Cauchy2} por Poisson em $1829$, \cite{Poisson} por Saint-Venant em $1837$, \cite{Saint-Venant1,Saint-Venant2,Saint-Venant3} e finalmente, de maneira mais convincente, por Stokes, em $1845$, \cite{Stokes1,Stokes2}. Enquanto o modelo sem viscosidade deduz-se de maneira muito simples de princípios físicos, o termo viscoso \'{e} muito mais delicado para justificar e, na aus\^{e}ncia de confirma\c{c}\~{a}o experimental s\'{o}lida, amplamente sujeito a controv\'{e}rsia. Foi ap\'{o}s Navier que as equa\c{c}\~{o}es de Navier-Stokes se firmaram como o modelo definitivo para o movimento de fluidos incompressíveis, desde escalas microsc\'{o}picas at\'{e} planet\'{a}rias.

Vamos primeiro formular precisamente o problema do Mil\^{e}nio associado \`{a}s equações de Navier-Stokes.  Primeiro, introduzimos o {\it espa\c{c}o de Schwartz} $\mathcal{S}(\mathbb{R}^3)$, o espa\c{c}o de fun\c{c}\~{o}es infinitamente diferenci\'{a}veis, tal que a fun\c{c}\~{a}o ela mesma, e todas as suas derivadas decaem mais r\'{a}pido do que qualquer polin\^{o}mio em infinito.

Fixe $\nu > 0$. Seja $u_0$ um campo vetorial de divergente nulo tal que suas componentes est\~{a}o em $\mathcal{S}(\mathbb{R}^3)$. Suponha tamb\'{e}m que o for\c{c}amento $F$ tenha componentes em $\mathcal{S}(\mathbb{R}^3 \times [0,\infty))$.

\vspace{0.3cm}

{\bf Problema do Mil\^{e}nio, \cite{MilPbmlNS}} - Prove que:

\vspace{0.5cm}

{\bf (N-S n\~{a}o forma singularidades em tempo finito)} Assuma que $F=0$. Mostre que, para qualquer $u_0 \in (\mathcal{S}(\real^3))^3$, existe uma solu\c{c}\~{a}o $u=u(x,t)$, cujas componentes pertencem a $C^{\infty}(\mathbb{R}^3 \times [0,+\infty))$, e uma press\~{a}o
$p \in C^{\infty}(\mathbb{R}^3 \times [0,+\infty))$, da equa\c{c}\~{a}o \eqref{EDPsFluidos} com velocidade inicial $u_0$, e tal que

\[\int_{\mathbb{R}^3} |u(x,t)|^2 \dd x \leq C < +\infty \text{ para todo } t \geq 0.\]

\vspace{0.3cm}

{\center \large \bf ou}

\vspace{0.3cm}

{\bf(N-S forma singularidades em tempo finito)} Ache  uma velocidade inicial $u_0  \in (\mathcal{S}(\real^3))^3$ e um for\c{c}amento $F \in (\mathcal{S}(\mathbb{R}^3 \times [0,\infty)))^3$
 para os quais não existe nenhuma solução  de \eqref{EDPsFluidos}, $u=u(x,t)\in (C^{\infty}(\mathbb{R}^3 \times [0,+\infty)))^3$, $p = p(x,t) \in C^{\infty}(\mathbb{R}^3 \times [0,+\infty))$, com este forçamento $F$ e este dado inicial $u_0$, e que satisfaz
\[\int_{\mathbb{R}^3} |u(x,t)|^2 \dd x \leq C < +\infty \text{ para todo } t \geq 0.\]

\vspace{0.3cm}

Para entender porque o problema se coloca, precisamos de contexto. Primeiro, observamos que temos existência de solu\c{c}\~{a}o para tempo curto.
Mais precisamente, suponha $u_0 \in H^m(\mathbb{R}^3)^3$, $\mbox{div }u_0 = 0$. Aqui, $H^m$ \'{e} \ o espa\c{c}o de fun\c{c}\~{o}es \ de quadrado integr\'{a}vel e cujas derivadas (no sentido fraco) at\'{e}  ordem $m$ s\~{a}o de quadrado integr\'{a}vel.

\vspace{0.3cm}

\begin{theorem}[Kato, 1972 \cite{Kato1972}] \label{ExistSuave}
Se $m > \frac{3}{2} + 2$ então, para algum $T>0$, existe uma \'{u}nica solu\c{c}\~{a}o
$u \in C([0,T];C^2(\mathbb{R}^3))^3 \cap C^1([0,T];C(\mathbb{R}^3))$ das equa\c{c}\~{o}es  de Navier-Stokes/Euler.
\end{theorem}

\vspace{.3cm}

Em dimens\~{a}o $N=2$ temos  $T = \infty$. 

O problema do Mil\^{e}nio tal como formulado acima tamb\'{e}m est\'{a} em aberto para as equa\c{c}\~{o}es de Euler, mas sua solu\c{c}\~{a}o, apesar de ser de grande interesse e relev\^{a}ncia, n\~{a}o tem um pr\^{e}mio associado. 

A quest\~{a}o em torno do problema do Mil\^{e}nio fica ainda mais interessante levando em conta o que acontece se abrirmos m\~{a}o de buscar solu\c{c}\~{o}es suaves, veja o resultado abaixo.

\vspace{.3cm}

\begin{theorem}[Leray, 1934 \cite{Leray1934} -- Hopf, 1951 \cite{Hopf1951}]
Suponha $u_0 \in (L^2(\mathbb{R}^3))^3$, $\dv u_0 = 0$ e $F=0$. Então existe (ao menos uma) {\em solução fraca} $u \in (L^\infty_{\loc} ((0,+\infty);L^2(\mathbb{R}^3)))^3 \cap (L^2_{\loc} ((0,+\infty);H^1(\mathbb{R}^3)))^3$ das equações de Navier-Stokes com dado inicial $u_0$. Para estas soluções vale
\[\frac{1}{2}\int_{\mathbb{R}^3}|u(x,T)|^2 \dd x + \nu \int_0^T \int_{\mathbb{R}^3} |Du (x,t)|^2 \dd x \dd t \leq \frac{1}{2} \int_{\mathbb{R}^3}|u_0(x)|^2 \dd x.\]
\end{theorem}

\vspace{.3cm}

Leray chamou estas soluções de ``soluções turbulentas".

\vspace{.3cm}

A unicidade de solu\c{c}\~{o}es fracas \'{e} um problema em aberto!

A press\~{a}o \'{e} tratada como um multiplicador de Lagrange, associada \`{a} incompressibilidade como vínculo.  Dada a velocidade, recupera-se a pressão tomando o divergente da equa\c{c}\~{a}o de velocidade e resolvendo o sistema elíptico resultante.

Para tornar clara a discuss\~{a}o acima, \'{e} necess\'{a}rio tornar preciso o que \'{e} uma solu\c{c}\~{a}o fraca.

\begin{definition}
Seja $u_0 \in (L^2(\mathbb{R}^3))^3$, $\dv u_0 = 0$ e $F=0$. Dizemos que $u \in (L^\infty_{\loc} ((0,+\infty);L^2(\mathbb{R}^3)))^3 \cap (L^2_{\loc} ((0,+\infty);H^1(\mathbb{R}^3)))^3$  é uma {\em solução fraca de Leray-Hopf} das equações de Navier-Stokes com dado inicial $u_0$ se:
\begin{enumerate}
  \item para todo $\Phi \in (C^\infty_c (\mathbb{R}^3 \times [0,+\infty)))^3$, $\dv \Phi (\cdot,t) = 0$, tivermos

  \[\int_0^{+\infty} \int_{\mathbb{R}^3} \partial_t \Phi \cdot u + [(u \cdot \nabla) \Phi] \cdot u + \int_{\mathbb{R}^3} \Phi(\cdot,0)\cdot u_0  = \nu \int_0^{+\infty} \int_{\mathbb{R}^3} \text{Tr}[(D\Phi) \, (Du)^t],\]

  \item $\dv u(\cdot,t) = 0$ no sentido das distribuições, para quase todo $t>0$,
  \item e se, para todo $t \geq 0$, valer a desigualdade de energia: \[\|u(\cdot,t)\|_{L^2}^2 \leq \|u_0\|_{L^2}^2 - 2\nu\int_0^t \|Du(\cdot,s)\|_{L^2}^2 \dd s.\]
\end{enumerate}
\end{definition}

Para a identidade integral acima fazer sentido, basta que a velocidade seja integr\'{a}vel no tempo, de quadrado integr\'{a}vel no espa\c{c}o.

Em $2011$ também foram obtidas soluções fracas globais para as equa\c{c}\~{o}es de Euler:

\begin{theorem}[Wiedemann, $2011$]
Suponha $u_0 \in (L^2(\torus^3))^3$, $\dv u_0 = 0$ e $F=0$. Então existem (infinitas!) soluções fracas $u \in (L^\infty_{\loc}((0,+\infty);L^2(\torus^3)))^3$ com este dado inicial. Mais ainda,
\[E(t) = \frac{1}{2} \int_{\torus^3} |u(x,t)|^2 \dd x \to 0 \text{ quando } t\to +\infty.\]
\end{theorem}

\vspace{.3cm}

As soluções obtidas por Wiedemann são muito exóticas: extremamente irregulares e seguem de uma construção bastante {\em ad hoc}, chamada integração convexa. Não são soluções fisicamente razoáveis.

Vamos explorar brevemente a natureza da dificuldade. A equação an\'{a}loga a Navier-Stokes em uma dimensão \'{e} a equa\c{c}\~{a}o de Burgers viscosa:

\[\partial_t u + u \partial_x u = \varepsilon \partial^2_{xx} u.\]

\vspace{.3cm}

Sabemos que, se $\varepsilon = 0$, então qualquer solu\c{c}\~{a}o com velocidade inicial suave, de suporte compacto que n\~{a}o seja identicamente nula {\em forma choques em tempo finito}. Mais precisamente, $\partial_x u \to +\infty$ em tempo finito.

\vspace{.3cm}

Contudo, se $\varepsilon > 0$ então estas soluç\~{o}es s\~{a}o globais no tempo. A possibilidade de forma\c{c}\~{a}o de singularidades est\'{a} conectada com um certo balan\c{c}o entre a n\~{a}o-linearidade e a difus\~{a}o, que varia conforme a dimens\~{a}o.

Para prosseguir, precisamos introduzir a {\it vorticidade} $\omega = \rote u$. Intuitivamente, a vorticidade \'{e} uma medida da velocidade de rota\c{c}\~{a}o infinitesimal do escoamento em cada ponto. Tomando-se o rotacional das equa\c{c}\~{o}es de Navier-Stokes e usando a condi\c{c}\~{a}o de divergente nulo deduz-se uma equa\c{c}\~{a}o de evolu\c{c}\~{a}o para a vorticidade:

\begin{equation}
\partial_t \omega + (u \cdot \nabla) \omega = (\omega \cdot \nabla ) u + \nu \Delta \omega.
\end{equation}
Esta equa\c{c}\~{a}o pode ser fechada acoplando-a com a condi\c{c}\~{a}o de divergente nulo e com a rela\c{c}\~{a}o $\omega = \rote u$.
O termo $\omega \cdot \nabla u$ pode ser reescrito como $Du \, \omega$, onde $Du$ \'{e} a derivada de $u$ em forma matricial.

Em dimens\~{a}o $2$, no caso $\nu = 0$, o lado direito se anula, pois $\omega$ \'{e} da forma $(0,0,\ast)$ e, na matriz $Du$, a terceira linha e terceira coluna se anulam. Consequentemente, a vorticidade se comporta como um escalar \'{e} \ {\it transportado} pelo escoamento. Mais precisamente, identificamos o vetor $\rote u$ com sua terceira componente, o escalar $\omega = \partial_{x_1}u_2 - \partial_{x_2} u_1$; nesse caso temos \[\partial_t \omega + u \cdot \nabla \omega = 0.\]

O transporte de vorticidade expresso pela equa\c{c}\~{a}o acima \'{e} suficiente para suprimir, em
dimens\~{a}o $2$, a forma\c{c}\~{a}o de singularidades, até mesmo se $\nu = 0$. Esse n\~{a}o \'{e} o caso em dimens\~{a}o $3$. O problema tem origem no termo $(\omega \cdot \nabla)u = Du \, \omega$ na equa\c{c}\~{a}o de vorticidade, que \'{e} comumente chamado de termo de {\it estiramento de vorticidade}. Para entender isso, note que $Du$ e $\omega$ s\~{a}o termos com a mesma ordem de regularidade. De fato, $Du$ pode ser expresso como um operador integral singular de tra\c{c}o nulo e ordem $0$ aplicado a $\omega$ (mais precisamente, $Du = D \nabla^{\perp} \Delta^{-1} \omega$). O operador $\partial_t + u\cdot\nabla$ \'{e} uma derivada temporal ao longo de trajet\'{o}rias de partículas, de modo que a equa\c{c}\~{o} de vorticidade, para $\nu = 0$, tem uma forma an\'{a}loga à de uma equa\c{c}\~{a}o de Ricatti, $\dot{W} = W^2$, que forma singularidades em tempo finito. Mesmo adicionando difus\~{a}o, a equa\c{c}\~{a}o de rea\c{c}\~{a}o-difus\~{a}o resultante ainda forma singularidades em tempo finito.

Uma complica\c{c}\~{a}o importante \'{e} a n\~{a}o-localidade do problema. O termo de estiramento $(\omega \cdot \nabla)u$ \'{e} quadr\'{a}tico em $\omega$, mas seu comportamento n\~{a}o pode ser reduzido ao que ocorre ao longo de uma trajet\'{o}ria de partícula.

Essa mesma n\~{a}o-localidade se reflete no tratamento da press\~{a}o na formula\c{c}\~{a}o do problema em termos de velocidade.

De fato, para evoluir a press\~{a}o \'{e} \ necess\'{a}rio resolver
\[-\Delta p = \dv \dv (u \otimes u),\]
para cada instante $t>0$.

\vspace{0.5cm}

Introduzimos o operador $K$, e escrevemos a lei de Biot-Savart: 
\[u = K[\omega], \qquad K = \mbox{ rot } (-\Delta)^{-1}.\]

Para explorar a relação do termo de estiramento de vorticidade com o problema de forma\c{c}\~{a}o de singularidades, \'{e} razoável
come\c{c}ar com problemas mais simples que apresentem características similares. Uma inst\^{a}ncia desta id\'{e}ia foi introduzir um modelo unidimensional para equa\c{c}\~{a}o de vorticidade, coloquialmente chamado {\it ``baby vorticity equation''}:

\begin{equation}
\partial_t \omega = \mathcal{H}(\omega) \omega,
\end{equation}
onde $\mathcal{H}(\omega)$ \'{e} \ a {\it transformada de Hilbert} de $\omega$,
\[\mathcal{H}(\omega) = \mathcal{H}(\omega)(x,t) = \frac{1}{\pi} VP \int \frac{\omega(y,t)}{x-y}\, dy,\]
com $x \in \mathbb{R}$, $t > 0$.

\vspace{0.5cm}

A integral $VP$ (valor principal) \'{e} um limite sim\'{e}trico em torno da singularidade em $x=y$. 

Para esta equação temos:

\begin{theorem} [Constantin-Lax-Majda, 1985 \cite{CLM1985}] Suponha $\omega_0=\omega(\cdot,0) \in H^1(\mathbb{R})$. Assuma tamb\'{e}m que
\[\{ x \;|\; \mathcal{H}(\omega_0) > 0 \} \neq \emptyset.\]
Seja $\omega=\omega(x,t)$ a solu\c{c}\~{a}o da ``baby vorticity equation''.
Ent\~{a}o existe $T>0$ tal que $\omega(x,t) \to \infty$ quando $t \to T$.
\end{theorem}

Usando propriedades elementares da transformada de Hilbert, \'{e} f\'{a}cil produzir exemplos de tais $\omega_0$.

Nos aproximando do problema do mil\^{e}nio, perguntamos o que acontece se adicionarmos um termo viscoso na ``baby vorticity''? De fato, as solu\c{c}\~{o}es desta vers\~{a}o viscosa tamb\'{e}m formam singularidades em tempo finito (Schochet, 1986 \cite{Schochet1986}).

\vspace{0.5cm}

Ainda buscando modelos com estrutura similar à das equa\c{c}\~{o}es de Euler e Navier-Stokes, mas agora em dimens\~{a} o $2$, consideramos din\^{a} mica de contorno; caso $\nu = 0$.

\vspace{0.5cm}

Considere um ``patch'' de vorticidade em dimens\~{a}o $2$, isto é, $\omega=\omega(x,t)=\chi_{\mathcal{D}(t)}$, onde $\mathcal{D}(t)$
 \'{e}  um {\it domínio material} (que se move com o escoamento).

\vspace{0.5cm}

Como a vorticidade, em dimensão $2$, é transportada pelo escoamento, a evolução do ``patch'' é descrita pela evolução da fronteira do ``patch''; seja $z = z(\alpha,t) \in \mathbb{R}^2$, $\alpha \in \mathbb{R}$, uma parametrização desta fronteira.

\vspace{0.5cm}

Então $z$ satisfaz uma equa\c{c}\~{a}o  de evolu\c{c}\~{a}o, chamda equação  do {\it contorno} do ``patch'', dada por
$\partial \mathcal{D}(t)$:
\begin{equation}
 \frac{dz}{dt}(\alpha,t) =
-\frac{1}{2\pi} \int_0^{2\pi} \log|z(\alpha,t)-z(\alpha',t)|  z_{\alpha}(\alpha',t)\,d \alpha'.
\end{equation}

Ap\'{o}s alguns c\'{a}lculos, é possível verificar que
\begin{equation}
 \frac{d z_{\alpha}}{dt}(\alpha,t) = \mathcal{V}(z_{\alpha}) z_{\alpha} ,
\end{equation}
onde
\[\mathcal{V}(z_{\alpha})= -\frac{1}{2\pi} VP \int_{\alpha - \pi}^{\alpha + \pi}z_{\alpha} (\alpha',t) \otimes \frac{z(\alpha,t)-z(\alpha',t)}{|z(\alpha,t)-z(\alpha',t)|^2} \,d \alpha'.\]

\vspace{0.5cm}

O operador $\mathcal{V}$ \'{e} um \ {\it operador integral singular}, de tra\c{c}o nulo e ordem $0$. Novamente, a equa\c{c}\~{a}o resultante tem uma estrutura an\'{a}loga \`{a} da equa\c{c}\~{a}o de vorticidade em dimens\~{a}o tr\^{e}s.

\vspace{0.5cm}

Coloca-se naturalmente a pergunta: formam-se singularidades no contorno em tempo finito?

\begin{theorem} [Chemin, 1993 \cite{Chemin1993}]
Assuma que $\partial \mathcal{D}_0 \in C^{1,\gamma}$, $0<\gamma<1$. Ent\~{a}o  $\partial \mathcal{D}(t)\in C^{1,\gamma}$, para todo $t\geq 0$.
Mais ainda, a curvatura e o comprimento do contorno s\~{a}o  estim\'{a}veis por uma {\em exponencial dupla} $\exp(\exp t)$.
\end{theorem}

\vspace{0.5cm}

Isso mostra que a vida n\~{a}o  \'{e}  t\~{a}o simples! Nem tudo que parece uma equa\c{c}\~{a}o de Ricatti, como a equa\c{c}\~{a}o de vorticidade 3D, forma singularidades.

Outro modelo amplamente estudado é a equa\c{c}\~{a}o SQG: equa\c{c}\~{a}o quase-geostr\'{o}fica de superficie, que emula o caso $\nu=0$:
\begin{equation} \label{2dqg}
\begin{array}{ll}
\partial_t \theta + u \cdot \nabla \theta = 0, & u = \nabla^{\perp}(-\Delta)^{-1/2} (-\theta).
\end{array}
\end{equation}
Acima, $\nabla^{\perp} = (-\partial_{x_2},\partial_{x_1})$.

\vspace{0.5cm}

Derivando, obtemos
\begin{equation}
\partial_t \,\nabla^{\perp}\theta + (u \cdot \nabla) \,\nabla^{\perp}\theta = D u \,\nabla^{\perp}\theta,
\end{equation}
logo $\nabla^{\perp}\theta$ faz o papel da vorticidade.

\vspace{0.5cm}

Caso dissipativo: adicione $\nu (-\Delta)^{\alpha}$ ao lado {\it esquerdo} de \eqref{2dqg}, $\alpha \geq 0$.
Caso critico: $\alpha = 1/2$. 

Kiselev {\it et alli} \cite{KNV2007} e Caffarelli {\it et alli} \cite{CV2010} mostraram, em trabalhos independentes, que n\~{a}o h\'{a} \ forma\c{c}\~{a}o de singularidades em tempo finito {\it para o caso critico} dissipativo. (Obs. O caso $\alpha > 1/2$ \'{e} mais fácil mostrar que não forma singularidade.)

\vspace{0.5cm}

O caso sem dissipa\c{c}\~{a}o (ainda) est\'{a} em aberto.

Mudando de assunto, vamos discutir brevemente a rela\c{c}\~{a}o entre mudan\c{c}as de escala e criticalidade.

Considere $u = u(x,t)$, $p=p(x,t)$, $x \in \mathbb{R}^N$, $t >0$, solução de Navier-Stokes.

\vspace{.3cm}

Seja $\lambda \in \mathbb{R}_+$.

\vspace{.3cm}

Então, \'{e} f\'{a}cil ver que
\[u_\lambda = u_\lambda (x,t) \equiv \lambda u(\lambda x, \lambda^2 t)\]
\[ p_\lambda = p_\lambda (x,t) \equiv  \lambda^2 p(\lambda x, \lambda^2 t)\]
também são soluções de Navier-Stokes.

\vspace{.3cm}

Essa invariância por mudan\c{c}a de escala sugere ``espaços críticos de dados iniciais".

Seja $u_0$ dado inicial. Então o ``scaling" de $u_0$ é:

\[(u_0)_\lambda = (u_0)_\lambda (x) = \lambda u_0 (\lambda x).\]

\vspace{.3cm}

Exemplo:
\begin{itemize}
  \item Se $N=2$ ent\~{a}o $u_0 \in L^2(\mathbb{R}^2) \Longrightarrow \|u_0\|_{L^2} = \|(u_0)_\lambda\|_{L^2}$ para todo $\lambda \in \mathbb{R}_+$.

  \item Se $N=3$ então $u_0 \in L^3(\mathbb{R}^3) \Longrightarrow \|u_0\|_{L^3} = \|(u_0)_\lambda\|_{L^3}$ para todo $\lambda \in \mathbb{R}_+$.
\end{itemize}

\vspace{.3cm}

Espaços críticos s\~{a}o aqueles cuja norma \'{e} invariante para as mudan\c{c}as de escala naturais do problema. Para Navier-Stokes, estes s\~{a}o: $L^2$ em dimens\~{a}o $N=2$ e $L^3$ em dimen\~{a}o $N=3$.

\vspace{.3cm}

Ao mesmo tempo, observamos a dificuldade adicional em dimens\~{a}o $N=3$ relativamente a $N=2$. O espaço $L^2$ \'{e} um {\em bom espa\c{c}o} -- é preservado pela evolu\c{c}\~{a}o -- enquanto $L^3$ {\em n\~{a}o \'{e}} preservado.

H\'{a} v\'{a}rios resultados de exist\^{e}ncia e unicidade local com dados iniciais em espa\c{c}os críticos.
Por exemplo, seja  $X_T$   o espaço dos campos com $\|\cdot\|_{X_T}$ finita, onde
\[ \|u\|_{X_T} \equiv \sup_{t\in (0,T)} t^{1/4}\|u(\cdot,t)\|_{L^6}.\]

\begin{theorem}[Kato 1984 \cite{Kato1984}]
Seja $u_0 \in L^3(\mathbb{R}^3)$ tal que $\dv u_0 = 0$. Então existe $T>0$ e uma \'{u}nica solu\c{c}\~{a}o $u \in (C(0,T;L^3(\mathbb{R}^3)))^3 \cap X_T$ das equa\c{c}\~{o}es de Navier-Stokes, com dado inicial $u_0$.

Mais ainda, existe $\varepsilon_0>0$ tal que, se $\|u_0\|_{L^3} < \varepsilon_0 \nu$, então podemos tomar $T=+\infty$.
\end{theorem}

Outros espa\c{c}os críticos com resultados semelhantes: $\dot{H}^{1/2}$ (Fujita Kato 1964 \cite{FK1964}), $BMO^{-1}$ (Koch Tataru 2001 \cite{KT2001}), sendo o espa\c{c}o $BMO^{-1}$ o maior espa\c{c}o crítico possível.

\vspace{.3cm}

Os teoremas de exist\^{e}ncia s\~{a}o, invariavelmente, {\em locais} para dados arbitr\'{a}rios, {\em globais} para dados pequenos. S\~{a}o resultados baseados em teoremas do ponto fixo via duas normas (m\'{e}todo de Kato). Vamos ser um pouco mais precisos.

\vspace{0.5cm}

Primeiro observamos que a ``solu\c{c}\~{a}o" no teorema de Kato \'{e} do tipo branda, ou ``mild", isto \'{e}, satisfaz uma forma integrada das equa\c{c}\~{o}es:

\[u = e^{\nu\Delta t}u_0 - \int_0^t e^{\nu(t-s)\Delta} \mathbb{P} [(u (\cdot,s) \cdot \nabla) u(\cdot,s)] \dd s. \]

Acima, $e^{\nu\Delta t}$ \'{e} o {\em semigrupo do calor}, de modo que, se $V= e^{\nu\Delta t}v_0 $ ent\~{a}o $V$ satisfaz a equa\c{c}\~{a}o do calor com dado inicial $v_0$:
\[\partial_t V = \nu\Delta V,\]
e $V(\cdot,0)= v_0$.

A no\c{c}\~{a}o de solu\c{c}\~{a}o mild \'{e} obtida aplicando o principio de Duhamel para resolver um problema do tipo
\[\partial_t V =  \nu\Delta V + F,\]
com $F = -\mathbb{P} [(u \cdot \nabla) u]$, e onde $\mathbb{P}$ \'{e} o {\it projetor de Leray}, a projeção do argumento nos campos de divergente nulo.

\vspace{0.5cm}

Note que  $\|u\|_{X_T} < + \infty$ significa
\[\|u(\cdot,t)\|_{L^6} \leq C t^{-1/4},\]
de modo que, eventualmente, $u(\cdot,t)$ se torna pequeno em $L^6$. Resultados parecidos com Kato $1984$ são até mais fáceis em $L^p$, $p>3$, de modo que, se $u(\cdot,t)$ se tornar {\em suficientemente pequeno} em $L^6$, seguir-se-á existência global.

\vspace{.3cm}

Este tipo de observação se aplica a vários outros espaços e nos leva a concluir que o problema de forma\c{c}\~{a}o de singularidades para Navier-Stokes \'{e} um problema de ``tempos intermedi\'{a}rios".

Uma linha de investiga\c{c}\~{a}o interessante \'{e} a busca por crit\'{e}rios para a forma\c{c}\~{a}o de singularidades, ou, condi\c{c}\~{o}es sobre as solu\c{c}\~{o}es que previnem sua forma\c{c}\~{a}o. O critério mais simples de formação de singularidades \'{e} devido a Serrin:

\begin{theorem}[Serrin 1962 \cite{Serrin1962}]
Seja $u_0 \in H^1(\mathbb{R}^3)$, $\mbox{div } u_0=0$. Seja $u \in L^\infty(0,T;H^1(\mathbb{R}^3))\cap L^2(0,T;H^2(\mathbb{R}^3))$ solu\c{c}\~{a}o de Navier-Stokes, onde $T<T_{\text{max}}$. Então, se $T_{\text{max}}<+\infty$ então
\[\lim_{T \to T_{\text{max}}^-} \;\; \int_0^T\|u(\cdot,t)\|_{L^q}^p \dd t \; = \; +\infty,\]
sempre que
\[\frac{2}{p}+\frac{3}{q} = 1, \; 2 \leq p \leq +\infty.\]
\end{theorem}

Apresentamos outro crit\'{e}rio conhecido como crit\'{e}rio BKM. Este \'{e} especialmente interessante por se aplicar tanto a Euler quanto a Navier-Stokes e ser adapt\'{a}vel a muitos outros problemas de interesse.

\begin{theorem}[Beale, Kato, Majda, 1984 \cite{BKM1984}]

\begin{enumerate}

\item Se $\int_0^T \| \omega(\cdot,s)\|_{L^{\infty}}\dd s$
for limitada, como fun\c{c}\~{a}o  de $T$, ent\~{a}o Euler (ou Navier-Stokes) n\~{a}o  forma
singularidade;

\item Se Euler (ou Navier-Stokes) formar singularidade em tempo $T_{\text{max}}$ ent\~{a}o
$\int_0^T \| \omega(\cdot,s)\|_{L^{\infty}} \dd s \rightarrow +\infty$ quando $T \nearrow T_{\text{max}}^-$.
\end{enumerate}
\end{theorem}

O crit\'{e}rio BKM \'{e} particularmente \'{u}til na investiga\c{c}\~{a}o num\'{e}rica de forma\c{c}\~{a}o de singularidades.

\vspace{0.5cm}

Ainda outro crit\'{e}rio -- a dire\c{c}\~{a}o  de vorticidade:  $\xi \equiv \omega / |\omega|$. Enquanto
\[\int_0^T \|D\xi (\cdot,s)\|_{L^{\infty}}^2 \dd s \]
for limitada n\~{a}o h\'{a} \ forma\c{c}\~{a}o \ de singularidades (Constantin, Fefferman, Majda, 1996 \cite{CFM1996}).

Note que os problemas de forma\c{c}\~{a}o de singularidades para Euler e Navier-Stokes est\~{a}o intimamente relacionados.
Em 2007 houve uma reuni\~{a}o cientifica para celebrar os $250$ anos da publica\c{c}\~{a}o \ das equa\c{c}~{o}es  de Euler
({\it Euler equations: 250 years on} \cite{Euler250}).

\vspace{0.5cm}

Houve uma sess\~{a}o de discuss\~{a}o sobre o Problema de Singularidades; alguns dias antes realizou-se uma pesquisa de opini\~{a}o. O resultado foi que a opini\~{a}o sobre forma\c{c}\~{a}o de singularidades para Euler era dividida mais-ou-menos ao meio, enquanto que, para Navier-Stokes, a maioria achava que n\~{a}o haveria forma\c{c}\~{a}o de singularidades. A situa\c{c}\~{a}o \'{e} muito diferente hoje.

Em outubro/2022 e maio/2023 Jiajie Chen e Tom Hou anunciaram uma demonstra\c{c}\~{a}o assistida por computador de forma\c{c}\~{a}o de singularidades para Euler.

\begin{theorem}[Chen Hou 2022, 2023 \cite{ChenHouPart1,ChenHouPart2}]
Existe uma familia de dados iniciais suaves, $(\theta_0,\omega_0)$, para os quais $2D$ Boussinesq e $3D$ Euler formam singularidades est\'{a}veis e quase auto-similares em tempo finito.
\end{theorem}

A assist\^{e}ncia do computador \'{e} necessária para construir os perfis auto-similares com erro pequeno e tamb\'{e}m para o c\'{a}lculo de majorantes \'{o}timos para as constantes no estudo de estabilidade.

O trabalho completo ainda n\~{a}o foi aceito para publica\c{c}\~{a}o.

\vspace{.3cm}

Em 2024 Hou anunciou que esta configuração sugere singularidade potencial também para $3D$ Navier-Stokes, \cite{Hou2024}. 

O ponto de partida deste trabalho \'{e} uma configura\c{c}\~{a}o axisim\'{e}trica e escoamento ``ímpar", isto \'{e}, a componente de rodopio da vorticidade tem simetria ímpar. A singularidade se origina na fronteira do domínio.

\vspace{1cm}

Para encerrar esta nota, enumeramos algumas conclus\~{o}es e conex\~{o}es

\begin{itemize}
 \item Os problemas de singularidades para Euler e Navier-Stokes são importantes e atuais.
 \item As ferramentas de an\'{a}lise usadas s\~{a}o modernas e sofisticadas: an\'{a}lise harm\^{o}nica (operadores integrais singulares,
decomposi\c{c}\~{a}o de Littlewood-Paley, ``wavelets'', etc), an\'{a}lise funcional (espa\c{c}os de fun\c{c}\~{o}es), an\'{a}lise microlocal
(c\'{a}lculo paradiferencial, etc).
 \item Esta \'{a}rea depende, de modo n\~{a}o-trivial, de sinergia entre an\'{a}lise, experimentos e simula\c{c}\~{o}es  num\'{e}ricas, na
 melhor tradi\c{c}\~{a}  de matem\'{a}tica aplicada.
 \item Uma resposta definitiva ajudar\'{a} \ a compreender v\'{a}rios outros problemas importantes de mec\^{a}nica dos fluidos, tais
como transi\c{c}\~{a}o \ para turbul\^{e}ncia, modelagem de escoamentos turbulentos e dissipa\c{c}\~{a}o de energia, etc.
 \item Em resumo, trata-se de uma \'{a}rea rica em problemas, interdisciplinar e extremamente ativa.
\end{itemize}



\bibliographystyle{plain}
\bibliography{bibMilPblm.bib}{}





\end{document}